\documentclass[11pt]{article}
\usepackage{amsmath}
\usepackage{amssymb}
\usepackage{amsthm}
\usepackage{amsfonts}
\usepackage{graphicx}
\usepackage{pdfpages}
\usepackage{xcolor,colortbl}
\usepackage{multicol}
\usepackage{url}
\usepackage{subfigure}

\newcommand\y{\cellcolor{green!10}}

\definecolor{Gray}{gray}{0.92}
\definecolor{LightCyan}{rgb}{0.88,1,1}

\newcolumntype{a}{>{\columncolor{Gray}}c}
\newcolumntype{b}{>{\columncolor{white}}c}

\usepackage{bbm,epsfig,graphics,epic,color,rotating,color}

    \usepackage[left=2.5cm,right=2.5cm,top=2.3cm,bottom=2.5cm]{geometry}

\numberwithin{equation}{section}
\def\E{{\mathbb E}}
\newcommand{\R}{{\mathbb R}}

\def\P{{\mathbb P}}
\def\N{{\mathbb N}}

\newtheorem*{theorem*}{Theorem}
\newtheorem{theorem}{Theorem}[section]

\newtheorem{corollary}[theorem]{Corollary}

\theoremstyle{definition}

\title{Multidimensional walks with random tendency}
\date{}
\author{Manuel Gonz\'alez-Navarrete}

\begin{document}

\maketitle %
\thispagestyle{empty} %
\baselineskip=14pt

\vspace{2pt}

\begin{abstract}

We introduce a multidimensional walk with memory and random tendency. The asymptotic behaviour is characterized, proving a law of large numbers and showing a phase transition from diffusive to superdiffusive regimes. In first case, we obtain a functional limit theorem to Gaussian vectors. In superdiffusive regime, we obtain strong convergence to a non-Gaussian random vector and characterize its moments.
\end{abstract}

%
\medskip
%


\section{Introduction}
\label{sec:intro}

The one-dimensional elephant random walk (ERW) was introduced in \cite{ST} (see also \cite{HK}). It can be represented by a sequence $\{X_1, X_2, \ldots\}$ where $X_i \in \{-1,+1\}$, for all $i\ge 1$. The position of the elephant at time $n$ is given by $S_{n}= \sum_{i=1}^{n}X_{i},$ and $S_0=0$. For the ERW model, it is supposed that the elephant remembers its full history and makes its $(n+1)$-step by choosing $t\in \{1,....n\}$ uniformly at random and then

\begin{equation}
\label{original}
X_{n+1} = \left\{
\begin{array}{ll}
X_{t}      &, \  \mbox{with probability $p$,} \\
-X_{t} &, \ \mbox{with probability $1-p$.}
\end{array}
\right.
\end{equation}
where $p \in [0,1]$ is a parameter.

In this sense, let denote $N(n,+1) = \#\{ i \in 1, 2, \ldots, n: X_i = +1\}$, the number of $+1$ steps until the $n$-step, and $N(n,-1)$ analogously defined. The position of the elephant can be rewritten as $S_n = N(n,+1) - N(n,-1)$. Therefore, the conditional probability of $(n+1)$-step in direction $+1$ is given by
\begin{eqnarray}
\label{Purn}
\P(X_{n+1} = 1 | N(n,+1), N(n,-1)) =  p \frac{N(n,+1)}{n} 
                                    + (1-p)\frac{N(n,-1)}{n},
\end{eqnarray}
{where $n=N(n,+1) + N(n,-1)$.} The probability \eqref{Purn} was used in \cite{BB} to introduce a relation with the classical P\'olya urn process \cite{Mah}. The evolution of a P\'olya urn is stated as follows. An urn starts with an initial quantity of $R_0$ red and $B_0$ blue balls and draws are made sequentially. After each draw, the ball is replaced and another ball of the same color is added to the urn. Let us consider the following notation: {the} urn is represented by the two-dimensional vector $(R_n,B_n)$, where $R_n$ and $B_n$ represent the number of red and blue balls at time $n\in \N$, respectively.

The relation established in \cite{BB} considers $N(n,+1) = R_n$, $N(n,-1) = B_n$ and the evolution in the P\'olya urn is modified in the sense that, at each draw the same color is used with probability $p$ and the opposite color with $1-p$. This fact allowed \cite{BB} to prove the existence of a transition from diffusive to super-diffusive behaviours for $S_n$, with critical $p_c=\frac{3}{4}$ (also proved by \cite{CGS}). That is, the mean squared displacement is a linear function of time in the diffusive case ($p < p_c$), but is given by a power law in the super-diffusive regime ($p > p_c$). Analogous results have been obtained in the $d$-dimensional formulation proposed in \cite{BL} (see also \cite{Ar,BerM,CVS,Mar}), by using martingale theory \cite{Ber,HH}. Such alternative approach was previously applied in the context of dependent Bernoulli sequences in \cite{BL,CGS,GL,WQY}.

In this paper we propose a generalization by introducing an evolution which employs the dependent dynamics from ERW and an independent one with random tendency. This kind of dynamics was introduced for one-dimensional models in \cite{GL} and also used in \cite{Ba,Bert,Bert2,Bu,KT}. In particular, at each time $n\ge 1$, the position is given by $S_n = (S_{n,1}, \ldots, S_{n,d})$ and the next step has two ways to be chosen. The election is given by a latent sequence $\{Y_n\}_{n\ge 1}$ of independent identically distributed Bernoulli random variables, with parameter $\theta$. We remark that the present formulation allows the possibility of delays (times for which the particle does not move), which has not been studied in dimensions $d\ge 2$. In addition, this process can be viewed as a model for the diffusion of $K$ opinions ($K = 2d$ or $2d+1$ with delays), generalizing \cite{GL1} and being an alternative to several formulations in the literature \cite{BPS,Ga,KHL,Mat}.

We characterize the asymptotic behaviours as function of parameters $p$ and $\theta$. The main results include a law of large numbers and functional limit theorems for diffusive and critical regimes, where $p_c = \frac{K+2\theta-1}{2\theta K}$. In the super-diffusive case, we stablish an almost sure convergence to a non-degenerated random vector, obtaining its first and second moments.

The rest of the paper is organized as follows. In Section \ref{sec:model} we introduce the walk and state the main results. Section \ref{sec:proof} includes the mathematical proofs based on generalized Polya urns (Theorems \ref{LGN} and \ref{continuous}) and martingale theory (Theorems \ref{LGNsuper} and \ref{moment}).

\section{Multidimensional random walks with memory and random tendency}
\label{sec:model}

We define a discrete-time evolution $(X_i)_{i \ge 1}$. The $n$-step denotes an opinion (movement), given by $X_n \in E= \{1,2, \ldots, K\}$ the set of choices. In the context of a random walk, we have $K = 2d$ or $2d+1$ with laziness, then, we denote the set of directions by

\begin{equation*}
E_d = \left\{
\begin{array}{ll}
(e_1, - e_1,  \ldots, e_d, - e_d)    &, \  \mbox{if K is even,} \\
(e_1, - e_1,  \ldots, e_d, - e_d, \mathbf{0})    &, \  \mbox{if K is odd,} 
\end{array}
\right.
\end{equation*}
where $(e_1, \ldots, e_d)$ is the canonical basis of the Euclidean space $\R^d$, and $\mathbf{0}$ denotes not movement.

Let $S_n = \displaystyle\sum_{i=1}^{n} X_i $ the $d$-dimensional position of the walker at time $n$. The $(n+1)$-step is obtained by flipping a coin with probability $\theta$, denoted $Y_n$ and then:

\begin{itemize}

\item If $Y_n=1$, we chose uniformly at random $t \in \{1,2, \ldots, n\}$, then $X_{n+1}$ is equal to $X_t$ with probability $p$. Otherwise, $X_{n+1}$ follows any other direction with uniform probability $\frac{1-p}{K-1}$.
\item If $Y_n=0$, then $X_{n+1}=e_1$ with probability $p$ or any other direction with uniform probability $\frac{1-p}{K-1}$.
\end{itemize}

Note that, if $\theta=1$ we obtain the dynamics defined by \cite{BL}. In case $\theta=0$, the tendency with intensity $p$ is given by direction $e_1$, such tendency is effective if $p > 1/K$.

Particular cases include: $K=2$, a one-dimensional walk with random tendency to the right (left) if $p> 1/2$ ($p<1/2$). Note that, if $\theta=0$ the central limit theorem guarantees that, for all $p \in [0,1]$, $\frac{1}{\sqrt{n}} (S_n - n(2p-1)) \to N(0,4p(1-p))$, as $n \to \infty$. For $\theta=1$, there exists a phase transition at $p_c=3/4$, from diffusive ($p < p_c$) to superdiffusive ($p > p_c$) behaviours, see \cite{BB}. In the case $K=3$, we have a random walk with delays or a three-opinion model (see for instance \cite{BPS,Ga,gut}). For larger $K$ the behaviours will be characterized in the following asymptotic results.

First, we obtain a strong law of large numbers for the mean position of the walker.


\begin{theorem}
\label{LGN}

Let $(S_n)_{n \in \mathbb{N}}$ the position of the walker, we get the following almost-surely convergence
\begin{equation*}
\lim_{n \to \infty} \frac{S_n}{n} = \frac{(1-\theta)(Kp-1)}{K-1 + \theta(1-K p)  } \left(1, 0, \ldots, 0\right)^T.
\end{equation*}
\end{theorem}

We also state the existence of a phase transition from diffusive to superdiffusive behaviours. We remark that the region of superdiffusivity is larger as $K$ increases. In diffusive and critical behaviours, we obtain a functional central limit theorem. We recall that this convergence holds on the function space $D[0, \infty )$ of right-continuous and left-bounded (c\`adl\`ag) functions, the Skorokhod space. Hereafter we will use the notation $\xrightarrow{d}$ for convergence in distribution (weak convergence), as $n$ diverges.

\begin{theorem} \label{continuous}
\label{teo2}
Let $(S_n)_{n \in \mathbb{N}}$, then

\begin{itemize}
\item[(i)] If $p < \frac{K+2\theta-1}{2\theta K}$ then, for $n \to \infty$, in $D[0, \infty)$

\begin{equation*}
\frac{1}{\sqrt{n}}\left[S_{\lfloor tn \rfloor} - \frac{tn(1-\theta)(Kp-1)}{K-1 + \theta(1-K p)  } \left(1, 0, \ldots, 0\right)^T \right] \xrightarrow{d} W_t,
\end{equation*}
where $W_t$ is a continuous $d$-dimensional Gaussian process with $W_0=(0,\ldots, 0)^T$, $\E(W_t)=(0,\ldots, 0)^T$ and, for $0 < s \le t$,

\begin{equation}
\label{corri}
\E(W_s W_t^T) = s\left(\frac{t}{s}\right)^{\frac{\theta(Kp-1)}{(K-1)}} \omega
\left( 
\begin{array}{cccccc}
(K+1) \alpha + \beta +p-1 & 0 &  & \ldots & & 0 \\
0 & 2 \beta & & \ldots & & 0 \\
 \vdots &  \vdots& & \ddots & & \vdots \\
0 & 0 & & \ldots & & 2\beta 
\end{array}
\right)
\end{equation}
where $ \omega=\dfrac{(K-1)(1-p)}{\beta^2(K-1 + 2\theta(1-K p) )}$, $\alpha = (K-1)p + \theta(1-K p)$ and $\beta = K-1 + \theta(1-K p).$

\medskip

\item[(ii)] If $p = \frac{K+2\theta-1}{2\theta K}$ then, for $n \to \infty$, in $D[0, \infty)$

\begin{equation*}
\frac{1}{\sqrt{n^t\log(n)}}\left[S_{\lfloor n^t \rfloor} - n^t \frac{K(2p-1)-1}{K -1}  \left(1, 0, \ldots, 0\right)^T \right] \xrightarrow{d} W_t,
\end{equation*}
where $W_t$ as above and for $0 < s \le t$,

\begin{equation}
\label{corrii}
\E(W_s W_t^T) =  4s\frac{1-p}{(K-1)^2}\left(p + \frac{K-3}{2} \right)\
\left( 
\begin{array}{cccccc}
(K+2)  & 0 &  & \ldots & & 0 \\
0 & 2 & & \ldots & & 0 \\
 \vdots &  \vdots& & \ddots & & \vdots \\
0 & 0 & & \ldots & & 2
\end{array}
\right)
\end{equation}
\end{itemize}
\end{theorem}

We remark that, letting $t=s=1$ it follow versions of the central limit theorem as in Theorems 3.3 and 3.6 in \cite{Ber}. In fact, Theorem \ref{continuous} generalizes the complementary results for multidimensional ERW in Theorems 4.2 and 4.3 from \cite{BerM}.

In superdiffusive case, the Gaussian behaviour is lost, but we have almost sure and mean square convergences to a non-degenerated random vector.

\begin{theorem}
\label{LGNsuper}

Let denote $\widehat{S}_n = S_n - \E(S_n)$ and $a=\frac{Kp-1}{K-1}$. If $p > \frac{K+2\theta-1}{2\theta K}$, then we have almost sure convergence

\begin{equation*}
\lim_{n \to \infty} \frac{\widehat{S}_n}{n^{a\theta}} = L,
\end{equation*}
where the limiting value $L$ is a non-degenerated random vector. We also have mean square convergence

\begin{equation*}
\lim_{n \to \infty} \E\left( \bigg\| \frac{\widehat{S}_n}{n^{a\theta}} - L \bigg\|^2 \right) = 0.
\end{equation*}
\end{theorem}

In addition, we obtain first and second moments of random vector $L$, which depend on the initial condition of the process.

\begin{theorem}
\label{moment}
The expected value of $L$ is $\E(L) = 0$, while its covariance matrix is obtained by
$$
\E(L L^T) = \displaystyle\lim_{n \to \infty} \frac{\Gamma(n)^2}{\Gamma(a\theta + n)^2} \E(\widehat{S}_n \widehat{S}_n^T),
$$
where
$$
\begin{array}{ll}
\E(\widehat{S}_n \widehat{S}_n^T) & = \displaystyle\prod_{i=1}^{n-1} \left( 1+ \frac{2a\theta}{i} \right) \E(\widehat{S}_1 \widehat{S}_1^T) + \displaystyle\sum_{i=1}^{n-2}\prod_{k=1}^{n-i} \left( 1+ \frac{2a\theta}{n+1-k} \right) \left[ \frac{\theta}{d} I_d + (1-\theta) M_p \right. \\
& \left. - \left(\displaystyle\frac{a\theta}{i}\prod_{l=1}^{i-1} \gamma_{i-l}\E(S_1) + (1-\theta)v_p\right)\left(\displaystyle\frac{a\theta}{i}\prod_{l=1}^{i-1} \gamma_{i-l}\E(S_1) + (1-\theta)v_p\right)^T \right]\\
&  + \frac{\theta}{d} I_d + (1-\theta) M_p + \displaystyle\prod_{k=1}^n \left( 1+ \frac{2a\theta}{n+1-k} \right) \E(\widehat{X}_1 \widehat{X}_1^T)\\
&  - \left(\displaystyle\frac{a\theta}{n-1}\prod_{l=1}^{n-2} \gamma_{n-1-l}\E(S_1) + (1-\theta)v_p\right)\left(\displaystyle\frac{a\theta}{n-1}\prod_{l=1}^{n-2} \gamma_{n-1-l}\E(S_1) + (1-\theta)v_p\right)^T.
\end{array}
$$
\end{theorem}

As a particular case, we include the following walk with uniform initial condition (see Theorem 3.8 of \cite{BL}).

\begin{corollary}
Let $\theta=1$ and $\E(S_1S_1^T)= \frac{1}{d} I_d$. Then, Theorem \ref{moment} implies $\E(L)=0$ and $\E(L L^T) = \frac{1}{d \Gamma(2a\theta + 1) \Gamma(2a\theta)} I_d$.
\end{corollary}

Finally, we remark that critical $p_c \ge \frac{1}{2}$, for all $\theta$ and $K$. Moreover, $p_c < 1$ if and only if, $\theta > \frac{1}{2}$. In particular, for the one-dimensional case see \cite{GL2} for details about critical curve, and \cite{Bert1} for further results in the diffusive regime.

\section{Proofs}
\label{sec:proof}

There are two elementary schemes to prove the phase transition and the limiting behaviours for this family of models. First, in one-dimensional case, the relation proposed in \cite{GL}, that is, the ERW behaves as a generalized binomial distribution \cite{DF}, then use classical results from martingale approach (see \cite{BL,Ber,CGS,HH,WQY}). Second, the theory from \cite{Jan} for the urn models with random replacement matrix, proposed by \cite{BB} and applied in \cite{GL1,GL2}.

On the one hand, the proofs of Theorems \ref{LGN} and \ref{continuous} are based on Janson \cite{Jan}. On the other hand, Theorems \ref{LGNsuper} and \ref{moment} are proved by using tools from Bercu and Laulin \cite{BL} based on martingale theory \cite{Du}.

\subsection{The strong law of large numbers and functional limit theorem}

Let denote $\mathcal{F}_n = \sigma(X_1, \ldots, X_n)$ the $\sigma$-field generated by the sequence $X_1, \ldots, X_n$. Then, for the $(n+1)$-step, choosing $t \in \{1, 2, \ldots, n\}$ uniformly, we use the following conditional probabilities for $x \in E_d$: 
\begin{equation}
\label{P1}
P(X_{n+1} = x | Y_n=1, \mathcal{F}_n) = \left\{
\begin{array}{ll}
p    &, \  \mbox{if $X_t= x$,} \\
\frac{1-p}{K-1} &, \ \mbox{otherwise,}
\end{array}
\right.
\end{equation}
and 

\begin{equation}
\label{P2}
P(X_{n+1} = x | Y_n=0, \mathcal{F}_n) = \left\{
\begin{array}{ll}
p    &, \  \mbox{if $x=e_1$,} \\
\frac{1-p}{K-1} &, \ \mbox{otherwise.}
\end{array}
\right.
\end{equation}

Therefore, defining $N(n,x) = | \{ i \in \{1,\ldots, n\}: X_i = x\}|$, the number of steps in the direction $x \in E_d$ until time $n$, we obtain 

\begin{equation}
\label{PFinal}
P(X_{n+1} = x | \mathcal{F}_n) = \left\{
\begin{array}{ll}
p   + \theta \left(\frac{1-Kp}{K-1}\right) \left( 1- \frac{N(n,e_1)}{n} \right)  &, \  \mbox{if $x=e_1$,} \\
\frac{1-p}{K-1}  + \theta \left(\frac{1-Kp}{K-1}\right) \frac{N(n,x)}{n}   &, \  \mbox{if $x\neq e_1$.}
\end{array}
\right.
\end{equation}

The position of the walker can be obtained by using an auxiliary process, which evolves as an urn model with $K$ colors. In this sense,

\begin{equation}
\label{relation}
S_n = \left\{
\begin{array}{ll}
(U_{1,n} - U_{2,n}, U_{3,n} - U_{4,n}, \ldots, U_{K-1,n}- U_{K,n})     &, \  \mbox{if $K$ is even}, \\[0.3cm]
(U_{1,n} - U_{2,n}, U_{3,n} - U_{4,n}, \ldots, U_{K-2,n}- U_{K-1,n})  &, \ \mbox{if $K$ is odd},
\end{array}
\right.
\end{equation}
where $U_n = (U_{1,n}, \ldots, U_{K,n})$ is the vector that denotes the number of balls of each of the $K$ colors, at time $n$. Each color is associated to the random variables $N(n,x)$ above.


Then, by defining the random replacement matrix in Section 2 of \cite{Jan}, we need to introduce the random vectors $\xi_i$, for $i \in \{1, \ldots, K\}$, which represent a random number of balls to be added into the urn. Essentially, these column vectors assume values on $\{e_1, \ldots, e_K\}$ the canonical basis of the Euclidean space $\R^K$. That is, these vectors denote the color of the ball to be added.

Therefore, using \eqref{PFinal} we obtain 
\begin{equation*}
\P(\xi_1 = e_1) = p \ \ \text{ and} \ \ \P(\xi_1 = e_i) = \frac{1-p}{K-1} \ \ \text{ for} \ i\neq 1.
\end{equation*}

Moreover, for $j \neq 1$,

$$
\P(\xi_j = x) = \left\{
\begin{array}{ll}
p+ \theta  \frac{1-Kp}{K-1}     &, \  \mbox{if } x=e_1, \\[0.2cm]
\frac{1-p - \theta(1-Kp)}{K-1} &, \  \mbox{if } x=e_j, \\[0.2cm]
 \frac{1-p}{K-1} &, \ \mbox{otherwise.}
\end{array}
\right.$$

In this sense, we obtain

\begin{equation}
\label{matrixA}
 A = \left( \E(\xi_1), \ldots, \E(\xi_K) \right) = \left(
\begin{array}{cccccc}
p & p + \theta\frac{1-Kp}{K-1} & & \ldots & &  p + \theta\frac{1-Kp}{K-1}  \\[0.3cm]
\frac{1-p}{K-1} & \frac{1-p - \theta(1-Kp)}{K-1} & & \ldots & & \frac{1-p}{K-1} \\[0.3cm]
\vdots & \vdots & & \ddots  & & \vdots \\[0.3cm]
\frac{1-p}{K-1} & \frac{1-p}{K-1} & & \ldots &  &\frac{1-p - \theta(1-Kp)}{K-1}
\end{array}
\right) \ ,
\end{equation}
for this matrix, the largest eigenvalue is $\lambda_1=1$, and for $j=2,\ldots,K$ we get
$$
\lambda_j = \theta \left(\frac{Kp-1}{K-1}\right).
$$

Now, we obtain vectors $u_1, \ldots, u_K$ and $v_1, \ldots, v_K$, such that

$$u_i^TA = \lambda_i u_i^T,$$
$$A v_i = \lambda_i v_i,$$
and
$$u_i^T v_j = \left\{
\begin{array}{ll}
1      &, \  \mbox{if } i=j, \\
0 &, \ \mbox{otherwise.}
\end{array}
\right.$$
That it, $u_1=(1,1,\cdots,1)^T$, and
$$v_1= ((K-1)(p - \lambda_2), 1-p, \ldots, 1-p)^T \frac{1}{(K-1)(1-\lambda_2)}.$$
Moreover, for $j=2,3,\ldots,K$ we obtain

$$
u_j=(1-p,\cdots,(K-1)\lambda_2-(K-2)-p,\cdots,1-p)^T \frac{1}{(K-1)(1-\lambda_2)},
$$
where the different value is at $j$-th position. Similarly, $v_j=(1,0,\ldots,-1,\ldots,0)^T,$ with $-1$ occupying the $j$-th position.

\textit{Proof Theorem \ref{LGN}}

We then use Theorem \textbf{3.21} from \cite{Jan}, which states that 
$$
n^{-1}U_n \longrightarrow \lambda_1 v_1 \ .
$$
Using \eqref{relation}, this finishes the proof.

\textit{Proof of Theorem \ref{continuous}}

For item $(i)$, we remark that Theorem \textbf{3.22} of \cite{Jan} says that the limiting covariance matrix is obtained as follows. First, let $L_I=\{i: \lambda_i < \lambda_1/2\}$ and $L_{II}=\{i: \lambda_i = \lambda_1/2\}$. From Collorary \textbf{5.3}(i) of \cite{Jan} we get

$$
\Sigma_I = \sum_{j,k \in L_I} \frac{u_j^TBu_k}{\lambda_1-\lambda_j-\lambda_k} v_jv_k^T  \ \ \ ; \ \ \ \Sigma_{II}=\sum_{j \in L_{II}} u_j^TBu_j v_jv_j^T,
$$

where $B=\displaystyle\sum_{i=1}^{K}v_{1i}B_i$ and $B_i=\E[\xi_i\xi_i^T]$, then $
B_1 = \left(
\begin{array}{cccc}
p&0&\cdots&0 \\
0&\frac{1-p}{K-1}&\cdots&0\\
\vdots& & \\
0&0&\cdots&\frac{1-p}{K-1}
\end{array}
\right),
$ and for $j\geq 2$ we get
$$
B_j = \left(
\begin{array}{cccccc}
p+\theta\left(\frac{1-Kp}{K-1}\right)&0&\cdots&\y{0}&\cdots&0 \\
0&\frac{1-p}{K-1}&\cdots&\y{0}&\cdots&0\\
\vdots&\vdots & \ddots &\y{\vdots}&\ddots & \vdots\\
\y{0}&\y{0}&\y{\cdots}&\y{\frac{1-p}{K-1}+\theta\left(\frac{Kp-1}{K-1}\right)}&\y{\cdots}&\y{0} \\
\vdots&\vdots &\ddots&\y{\vdots}&\ddots& \vdots\\
0&0&\cdots&\y{0}&\cdots&\frac{1-p}{K-1}
\end{array}
\right),
$$
where here we have highlighted $j$-column and $j$-row. A direct computation leads to
$$
B=\frac{1}{(K-1)+ \theta(1-Kp)} \left(
\begin{array}{cccc}
p(K-1)+\theta(1-Kp)&0&\cdots&0 \\
0&1-p&\cdots&0 \\
\vdots&\vdots & \ddots& \\
0 & 0 & \cdots & 1-p
\end{array}
\right),
$$

Therefore,

$$u_i^T  B u_j =  \frac{1-p}{(K-1)^2 (1-\lambda_2)^2} \cdot \left\{
\begin{array}{ll}
 p-1     &, \  \mbox{if } i\neq j, \\
 p-1 + (K-1)(1-\lambda_2) &, \ \mbox{if } i= j,
\end{array}
\right.$$
and
$$
v_i v_j^T = \left(
\begin{array}{cccccc}
1&0&\cdots&\y{-1}&\cdots&0 \\
0&0&\cdots&\y{0}&\cdots&0\\
\vdots&\vdots & \ddots &\y{\vdots}&\ddots & \vdots\\
\y{-1}&\y{0}&\y{\cdots}&\y{1}&\y{\cdots}&\y{0} \\
\vdots&\vdots &\ddots&\y{\vdots}&\ddots& \vdots\\
0&0&\cdots&\y{0}&\cdots&0
\end{array}
\right),
$$
we highlighted $j$-column and $i$-row.
Finally, we obtain

\begin{equation}
\label{SigmaI}
\Sigma_I= C \left(
\begin{array}{cccccccc}
(K-1)a& -a & & -a & & \ldots & &  -a\\[0.3cm]
-a & b & &p-1 & & \ldots & & p-1 \\[0.3cm]
-a  &p-1& & b& & \ldots & & p-1 \\[0.3cm]
\vdots & \vdots & & \vdots & & \ddots  & & \vdots \\[0.3cm]
-a & p-1 & & p-1 & & \ldots &  & b
\end{array}
\right) \ ,
\end{equation}

where $a=(K-1) (p-\lambda_2), \ b=(p-1)+(K-1)(1-\lambda_2)$ and $C=\frac{1-p}{(K-1)^2 (1-\lambda_2)^2(1-2\lambda_2)}$. This is the limiting covariance matrix for the proportions $U_n$. Then, using \eqref{relation}, we finish the proof.

Item $(ii)$. For this we use Theorem \textbf{3.23} and Corollary \textbf{5.3-(i)} of \cite{Jan} to obtain the following limiting covariance matrix to finish the proof
\begin{equation}
\label{SigmaII}
\Sigma_{II}= 4\frac{1-p}{(K-1)^2}\left(p + \frac{K-3}{2} \right)\left(
\begin{array}{cccccccc}
(K-1)& -1 & & -1 & & \ldots & &  -1\\[0.3cm]
-1 & 1 & &0 & & \ldots & & 0 \\[0.3cm]
-1  &0 & & 1& & \ldots  & &0 \\[0.3cm]
\vdots & \vdots & & \vdots & & \ddots  & & \vdots \\[0.3cm]
-1 & 0 & & 0 & & \ldots &  & 1
\end{array}
\right) \ .
\end{equation}

\subsection{Superdiffusive behaviours}

We define a locally square-integrable multidimensional martingale, given by

\begin{equation}
M_n = a_n \widehat{S}_n = \displaystyle\sum_{k=1}^{n} a_k \left(\widehat{S}_k - \left(1+\frac{a\theta}{k-1}\right)\widehat{S}_{k-1}\right) = \displaystyle\sum_{k=1}^{n} a_k\varepsilon_k,
\end{equation}
where $\widehat{S}_n = S_n - \E(S_n)$, $a_k = \displaystyle\prod_{l=1}^{k-1} \frac{l}{l+a\theta}$ and $a=\frac{Kp-1}{K-1}$.

\textit{Proof of Theorem \ref{LGNsuper}.}

Then, by Theorem \textbf{4.3.15} from \cite{Du}, we need to prove that

\begin{equation}
\displaystyle\lim_{n \to \infty} Tr \langle M \rangle_n < \infty \ \ \text{ a.s}
\end{equation}
where $Tr A$ stands for the trace of matrix $A$ and

\begin{equation}
\langle M \rangle_n = \displaystyle\sum_{k=1}^{n} \E \big[ (a_k \varepsilon_k)(a_k \varepsilon_k)^T \big| \mathcal{F}_{k-1}\big].
\end{equation}
In this sense, we obtain

$$
\E(\varepsilon_{n+1}\varepsilon_{n+1}^T \big| \mathcal{F}_n )  =  -\left(\frac{a\theta}{n}\right)^2 \widehat{S}_n \widehat{S}_n^T + \E( \widehat{X}_{n+1}  \widehat{X}_{n+1}^T\big| \mathcal{F}_n ),$$

where

\begin{equation}
\label{XX}
\begin{array}{lll}
\E(\widehat{X}_{n+1}\widehat{X}_{n+1}^T \big| \mathcal{F}_n )  & = & \E({X}_{n+1}{X}_{n+1}^T \big| \mathcal{F}_n ) - \frac{a\theta}{n} \left(\frac{a\theta}{n} S_n + (1-\theta)v_p \right)\left( \E(S_n) - S_n \right)^T \\[0.3cm]
&& - \left(\frac{a\theta}{n} S_n + (1-\theta)v_p \right)\left(\frac{a\theta}{n} \E(S_n) + (1-\theta)v_p \right)^T.
\end{array}
\end{equation}
Finally,
\begin{equation}
\begin{array}{lll}
\E(\varepsilon_{n+1}\varepsilon_{n+1}^T \big| \mathcal{F}_n )  & = & (1-\theta) M_p + \frac{a\theta}{n}\Sigma_n + \frac{\theta(1-a)}{d} I_d \\[0.3cm]
&& - \left(\frac{a\theta}{n} S_n + (1-\theta)v_p \right)\left(\frac{a\theta}{n} S_n + (1-\theta)v_p \right)^T,
\end{array}
\end{equation}
where $
M_p= \left(
\begin{array}{cccc}
p& 0&  \ldots & 0\\[0.3cm]
 0& \frac{1-p}{K-1}&  \ldots & 0\\[0.3cm]
 \vdots &  \vdots &  \ddots   & \vdots \\[0.3cm]
0& 0&  \ldots & \frac{1-p}{K-1}
\end{array}
\right)$
and $\Sigma_n = \displaystyle\sum_{i=1}^{d} N^X_n(i) e_i e_i^T$, where $ N^X_n(i)=\displaystyle\sum_{k=1}^n  \mathbb{I}_{X_k^i \neq 0}$ with $X_n^i$ being the $i$-th coordinate of $X_n$, which denotes the number of steps in such direction until time $n$. Moreover, $v_p = \left(p, \frac{1-p}{K-1}, \ldots, \frac{1-p}{K-1}\right)^T$.

Then,
\begin{equation}
\begin{array}{lll}
Tr \langle M \rangle_n & = & a_1^2 \E(\varepsilon_1 \varepsilon_1^T) + \alpha(p,\theta,K) \displaystyle\sum_{l=1}^{n-1} a_{l+1}^2 \left(1- \frac{2\theta(1-\theta)}{\alpha(p,\theta,K)}  \frac{Tr (S_l v_p^T)}{l} \right) \\[0.3cm]
&& - a^2\theta^2\displaystyle\sum_{l=1}^{n-1} \left(\frac{a_{l+1}}{l}\right)^2 \| S_l \|^2,
\end{array}
\end{equation}
where $\alpha(p,\theta,K) = 1-(1-\theta)^2 \left( p \frac{(1-p)^2}{K-1}\right)$. In addition, note that for all $l \ge 1$ and for all $p, K$ and $\theta$, $-1 \le \frac{2\theta(1-\theta)}{\alpha(p,\theta,K)}  \frac{Tr (S_l v_p^T)}{l} \le 1$. Then,

\begin{equation}
Tr \langle M \rangle_n  \le   a_1^2 \E(\varepsilon_1 \varepsilon_1^T) + 2 \alpha(p,\theta,K) \displaystyle\sum_{l=1}^{n-1} a_{l+1}^2.
\end{equation}

Note that, $ \displaystyle\sum_{l=1}^{n} a_{l}^2 =\displaystyle\sum_{l=1}^{n} \left(\frac{\Gamma(a\theta+1) \Gamma(l)}{\Gamma(a\theta +l)}\right)^2$, which in the superdiffusive regime satisfies

\begin{equation}
\displaystyle\lim_{n \to \infty} \displaystyle\sum_{l=1}^{n} a_{l}^2 = \displaystyle\sum_{l=1}^{\infty} \left(\frac{\Gamma(a\theta+1) \Gamma(l)}{\Gamma(a\theta +l)}\right)^2 = {}_3F_2\left(1,1,1;a\theta+1,a\theta+1;1\right),
\end{equation}
the finite confluent hypergeometric function. Therefore,

\begin{equation}
\displaystyle\lim_{n \to \infty} Tr \langle M \rangle_n  <   \infty \ \text{ a.s }.
\end{equation}

Finally, using $L_n = \frac{M_n}{\Gamma(a\theta + 1)}$, second part of Theorem \textbf{4.3.15} of \cite{Du} guarantees

\begin{equation}
\displaystyle\lim_{n \to \infty} M_n = M = \displaystyle\sum_{l=1}^{\infty} a_{l}\varepsilon_l \text{ and } \displaystyle\lim_{n \to \infty} L_n = L, 
\end{equation}
almost surely. In addition with property $\displaystyle\lim_{n \to \infty}  \frac{\Gamma(n + a\theta)}{\Gamma(n)n^{a\theta}} =1$, then $\displaystyle\lim_{n \to \infty}  n^{a\theta} a_n = \Gamma(a\theta + 1)$. Therefore, $\frac{\widehat{S}_n}{n^{a\theta}} = \frac{\Gamma(1 + a\theta) L_n}{a_n n^{a\theta}} \to L$ almost surely, as $n \to \infty$.

The mean square convergence follows since
\begin{equation}
\E( \|M_n\|^2 ) = \E( Tr \langle M \rangle_n) \le c  \displaystyle\sum_{l=1}^{n} a_l^2.
\end{equation}

Hence, $\displaystyle\sup_{n \ge 1} \E(\|M_n\|^2 ) \le c \cdot {}_3F_2\left(1,1,1;a\theta+1,a\theta+1;1\right) < \infty$. Which means $(M_n)$ is a bounded martingale in $\mathbb{L}^2$. Then, the result holds.

\textit{Proof of Theorem \ref{moment}.}

Note that $\E(M_n)=0$ for all $n \ge 1$, then $\E(M)=0$ and $\E(L)=0$. We also know that

$$
\begin{array}{lll}
\E(\widehat{S}_{n+1} \widehat{S}_{n+1}^T) & = & \left( 1+ \frac{2a\theta}{n} \right) \E(\widehat{S}_n \widehat{S}_n^T) + \E(\widehat{X}_{n+1} \widehat{X}_{n+1}^T)\\[0.3cm]
&  =  &\displaystyle\prod_{i=1}^n \left( 1+ \frac{2a\theta}{i} \right) \E(\widehat{S}_1 \widehat{S}_1^T)   +  \E(\widehat{X}_{n+1}\widehat{X}_{n+1}^T)\\[0.4cm]
& & +\displaystyle\sum_{i=0}^{n-1} \prod_{k=1}^{n-i}\left(1+\displaystyle\frac{2a\theta}{n+1-k}\right) \E(\widehat{X}_{i+1}\widehat{X}_{i+1}^T),
\end{array}
$$
by taking expectation in both sides of \eqref{XX}, using $\E(S_n) = \displaystyle\prod_{i=1}^{n-1}\gamma_{n-i} S_1$, and $\E(\Sigma_n) = \frac{n}{d} I_d$ from (5.24) in \cite{BL}, then we finish the proof.

\subsection*{Acknowledgements}
The author thanks Rodrigo Lambert and Eugene Pechersky for several comments. This work was partially supported by Fondo Especial DIUBB 1901083-RS from Universidad del B\'io-B\'io.

\medskip

{\scriptsize
{\sc Departamento de Estad\'i{}stica, Universidad del B\'io-B\'io. Avda. Collao 1202, CP 4051381, Concepci\'on, Chile. E-mail address: magonzalez@ubiobio.cl
}

\end{document}